
%
%
%
%
%
%
\magnification=\magstephalf      
%
%
\vsize=7.5truein                 
\hsize=5.2truein                 
\newskip\stdskip                 
\stdskip=6pt plus3pt minus3pt    
\medskipamount=\stdskip          
\parindent=0pt                   
\parskip=\stdskip                
\abovedisplayskip=\stdskip       
\belowdisplayskip=\stdskip       
\mathsurround=0.75pt             
\overfullrule=0pt                
%
%
\def\ppar{\par\goodbreak\vskip 8pt plus 4pt minus 4pt}     
%
%
\def\stdspace{\hskip 0.75em plus 0.15em\ignorespaces}
\let\qua\stdspace 
%
%
%
%
%
%
%
\def\hexnumber#1{\ifcase#1 0\or 1\or 2\or 3\or 4\or 5\or 6\or 7\or 8\or
 9\or A\or B\or C\or D\or E\or F\fi}
%
%
\font\thirtnmsa=msam10 scaled 1315    
\font\tenmsa=msam10          \font\ninemsa=msam9
\font\sevenmsa=msam7         \font\sixmsa=msam6
\font\fivemsa=msam5
%
%
\newfam\msafam                  \textfont\msafam=\tenmsa
\scriptfont\msafam=\sevenmsa    \scriptscriptfont\msafam=\fivemsa
\edef\hexa{\hexnumber\msafam}        
\def\msa{\fam\msafam\tenmsa}         
%
%
\font\thirtnmsb=msbm10 scaled 1315   
\font\tenmsb=msbm10      \font\ninemsb=msbm9
\font\sevenmsb=msbm7     \font\sixmsb=msbm6
\font\fivemsb=msbm5
%
\newfam\msbfam                   \textfont\msbfam=\tenmsb       
\scriptfont\msbfam=\sevenmsb     \scriptscriptfont\msbfam=\fivemsb
\edef\hexb{\hexnumber\msbfam}    
\def\msb{\fam\msbfam\tenmsb}     
%
%
\font\thirtneufm=eufm10 scaled 1315   
\font\teneufm=eufm10                 \font\nineeufm=eufm9
\font\seveneufm=eufm7                \font\sixeufm=eufm6
\font\fiveeufm=eufm5
%
\newfam\eufmfam                    \textfont\eufmfam=\teneufm
\scriptfont\eufmfam=\seveneufm     \scriptscriptfont\eufmfam=\fiveeufm
\edef\hexf{\hexnumber\eufmfam}      
\def\frak{\fam\eufmfam\teneufm}     
%
%
%
\font\thirtnrm=cmr10 scaled 1315    
\font\ninerm=cmr9                   \font\sixrm=cmr6   
%
\font\thirtni=cmmi10 scaled 1315    
\font\ninei=cmmi9                   \font\sixi=cmmi6  
%
\font\thirtnsy=cmsy10 scaled 1315   
\font\ninesy=cmsy9                  \font\sixsy=cmsy6  
%
\font\thirtnbf=cmbx10 scaled 1315   
\font\ninebf=cmbx9                  \font\sixbf=cmbx6  
%
%
\font\thirtnex=cmex10 scaled 1315   
\font\nineex=cmex9                  
%
%
\font\thirtnit=cmti10 scaled 1315  
\font\nineit=cmti9                  
%
\font\thirtnsl=cmsl10 scaled 1315  
\font\ninesl=cmsl9                  
%
\font\thirtntt=cmtt10 scaled 1315  
\font\ninett=cmtt9                  
%
%
%
%
\def\small{%
%
%
\textfont0=\ninerm \scriptfont0=\sixrm \scriptscriptfont0=\fiverm
\def\rm{\fam0\ninerm}
%
%
\textfont1=\ninei \scriptfont1=\sixi \scriptscriptfont1=\fivei
%
%
\textfont2=\ninesy \scriptfont2=\sixsy \scriptscriptfont2=\fivesy
%
%
\textfont3=\nineex \scriptfont3=\nineex \scriptscriptfont3=\nineex
%
%
\textfont\bffam=\ninebf \scriptfont\bffam=\sixbf
\scriptscriptfont\bffam=\fivebf \def\bf{\fam\bffam\ninebf}%
%
%
\textfont\itfam=\nineit \def\it{\fam\itfam\nineit}%
\textfont\slfam=\ninesl \def\sl{\fam\slfam\ninesl}%
\textfont\ttfam=\ninett \def\tt{\fam\ttfam\ninett}%
%
%
%
\textfont\msafam=\ninemsa \scriptfont\msafam=\sixmsa
\scriptscriptfont\msafam=\fivemsa \def\msa{\fam\msafam\ninemsa}%
%
%
\textfont\msbfam=\ninemsb \scriptfont\msbfam=\sixmsb
\scriptscriptfont\msbfam=\fivemsb \def\msb{\fam\msbfam\ninemsb}%
%
%
\textfont\eufmfam=\nineeufm  \scriptfont\eufmfam=\sixeufm
\scriptscriptfont\eufmfam=\fiveeufm \def\frak{\fam\eufmfam\nineeufm}%
%
%
%
\normalbaselineskip=11pt%
\setbox\strutbox=\hbox{\vrule height8pt depth3pt width0pt}%
%
%
\normalbaselines\rm
%
%
\stdskip=4pt plus2pt minus2pt    
\medskipamount=\stdskip          
\parskip=\stdskip                
\abovedisplayskip=\stdskip       
\belowdisplayskip=\stdskip       
\def\ppar{\par\goodbreak\vskip 6pt plus 3pt minus 3pt}%
%
%
\def\section##1{\global\advance\sectionnumber by 1
\vskip-\lastskip\penalty-800\vskip 20pt plus10pt minus5pt 
\egroup{\bf\number\sectionnumber\quad##1}\bgroup\small         
\vskip 6pt plus3pt minus3pt
\nobreak\resultnumber=1}
}    
%
\def\beginsmall{\bgroup\small}
\let\endsmall\egroup
%
%
%
%
\def\large{%
\textfont0=\thirtnrm \scriptfont0=\ninerm \scriptscriptfont0=\sevenrm
\def\rm{\fam0\thirtnrm}%
\textfont1=\thirtni \scriptfont1=\ninei \scriptscriptfont1=\seveni
\textfont2=\thirtnsy \scriptfont2=\ninesy \scriptscriptfont2=\sevensy
\textfont3=\thirtnex \scriptfont3=\thirtnex \scriptscriptfont3=\thirtnex
\textfont\bffam=\thirtnbf \scriptfont\bffam=\ninebf
\scriptscriptfont\bffam=\sevenbf \def\bf{\fam\bffam\thirtnbf}%
\textfont\itfam=\thirtnit \def\it{\fam\itfam\thirtnit}%
\textfont\slfam=\thirtnsl \def\sl{\fam\slfam\thirtnsl}%
\textfont\ttfam=\thirtntt \def\tt{\fam\ttfam\thirtntt}%
\textfont\msafam=\thirtnmsa \scriptfont\msafam=\ninemsa
\scriptscriptfont\msafam=\sevenmsa \def\msa{\fam\msafam\thirtnmsa}%
\textfont\msbfam=\thirtnmsb \scriptfont\msbfam=\ninemsb
\scriptscriptfont\msbfam=\sevenmsb \def\msb{\fam\msbfam\thirtnmsb}%
\textfont\eufmfam=\thirtneufm  \scriptfont\eufmfam=\nineeufm
\scriptscriptfont\eufmfam=\seveneufm \def\frak{\fam\eufmfam\teneufm}%
\normalbaselineskip=16pt%
\setbox\strutbox=\hbox{\vrule height11.5pt depth4.5pt width0pt}%
\normalbaselines\rm}%
\let\Large\large   
%
\def\Bbb#1{{\msb#1}}

%

\def\re{\Bbb R}
%
\mathchardef\plussquare="0\hexa01
\mathchardef\nge="3\hexb0B
\mathchardef\maltesecross="0\hexa7A
\mathchardef\del="0\hexf01
%
%
%
%
\font\sc=cmcsc10
%
%
%
%
\def\sqr#1#2{{\vcenter{\vbox{\hrule  height.#2truept
	\hbox{\vrule width.#2truept height#1truept 
	\kern#1truept \vrule width.#2truept}
	\hrule height.#2truept}}}}
\def\sq{\sqr55}    
%
%
%
%
\newcount\sectionnumber            
\newcount\resultnumber             
\sectionnumber=0\resultnumber=1    
%
%
%
\def\section#1{\global\advance\sectionnumber by 1
\xdef\nextkey{\number\sectionnumber}
\vskip-\lastskip\penalty-800\vskip 20pt plus10pt minus5pt 
{\large\bf\number\sectionnumber\quad#1}         
\vskip 8pt plus4pt minus4pt
\nobreak\resultnumber=1}      
%
%
%
%
%
\def\sh#1{\vskip-\lastskip\ppar{\bf #1}\par\nobreak\medskip}         
%
%
%
%

%
\def\proc#1{\xdef\nextkey{\number\sectionnumber.\number\resultnumber}%
\vskip-\lastskip\ppar\bf%
\noindent#1\ \number\sectionnumber.\number\resultnumber
\stdspace\sl\global\advance\resultnumber by 1\ignorespaces}
\def\endproc{\rm\ppar} 
%
%
\def\prf{\vskip-\lastskip\ppar\noindent{\bf Proof}%
\stdspace\rm}                            
\def\qed{\hfill$\sq$\par\goodbreak\rm}   
\def\endprf{\unskip\stdspace\hbox{}
\hfill$\sq$\par\medskip}                 
%
%
%
%
%
%
%
%
\def\proclaim#1{\vskip-\lastskip\ppar\bf%
\noindent#1\stdspace\sl\ignorespaces} 

%
%
%
%
\def\rk#1{\vskip-\lastskip\ppar{\bf #1}\stdspace\ignorespaces}                

%
%
%
%
%
%
\def\label{\xdef\nextkey{\number\sectionnumber.\number\resultnumber}%
\number\sectionnumber.\number\resultnumber
\global\advance\resultnumber by 1}
%
%
%
%
%
%
%
%
%
%
%
%
%
%
%
%
\newcount\refnumber              
\refnumber=1                     
\long\def\reflist#1\endreflist{%
\long\def\thereflist{#1}{\def\refkey##1##2\par{\xdef##1{\number\refnumber}%
\global\advance\refnumber by 1}%
\def\key##1##2\par{\expandafter\xdef%
\csname##1\endcsname{\number\refnumber}%
\global\advance\refnumber by 1}#1\par}}
\long\def\references{%
\penalty-800\vskip-\lastskip\vskip 15pt plus10pt minus5pt 
{\large\bf References}\ppar 
{\leftskip=25pt\frenchspacing    
\small\parskip=3pt plus2pt       
\def\refkey##1##2\par{\noindent  
\llap{[##1]\stdspace}\ignorespaces##2\par}         
\def\key##1##2\par{\noindent  
\llap{[\ref{##1}]\stdspace}\ignorespaces##2\par}  
\def\,{\thinspace}\thereflist\par}}
%
%
%
\newcount\footnotenumber         
\footnotenumber=1                
\def\fnote#1{\xdef\nextkey{\number\footnotenumber}%
{\small\ifnum\footnotenumber>9\parindent=14pt%
\else\parindent=10pt\fi\footnote{$^{\number\footnotenumber}$}%
{\hglue-5pt#1}\global\advance\footnotenumber by 1}}
%
%
%
%
%
%
%
\newcount\figurenumber          
\figurenumber=1                 
\def\caption#1{\xdef\nextkey{\number\figurenumber}%
\cl{\small Figure \number\figurenumber: #1}%
\global\advance\figurenumber by 1}
\def\figurelabel{\xdef\nextkey{\number\figurenumber}%
\cl{\small Figure \number\figurenumber}%
\global\advance\figurenumber by 1}
\long\def\figure#1\endfigure{{\xdef\nextkey{\number\figurenumber}%
\let\captiontext\relax\def\caption##1{\xdef\captiontext{##1}}%
\midinsert\cl{\ignorespaces#1\unskip\unskip\unskip\unskip}\vglue6pt\cl{\small 
Figure \number\figurenumber\ifx\captiontext\relax\else: \captiontext
\fi}\endinsert\global\advance\figurenumber by 1}}
%
%
%
%
%
%
%
\def\nextkey{??}   
%
\def\key#1{\expandafter\xdef\csname #1\endcsname{\nextkey}}
\def\ref#1{\expandafter\ifx\csname #1\endcsname\relax
\immediate\write16{Reference {#1} undefined}??\else
\csname #1\endcsname\fi}
%
%
%
%
%
%
%
\newread\gtinfile
\newwrite\gtreffile
\def\useforwardrefs{
\openin\gtinfile\jobname.ref
\ifeof\gtinfile
\closein\gtinfile
\immediate\write16{No file \jobname.ref}
\else
\closein\gtinfile
\input \jobname.ref
\fi
\immediate\openout\gtreffile \jobname.ref
%
%
\def\key##1{{\def\\{\noexpand}%
\expandafter\xdef\csname ##1\endcsname{\nextkey}%
\immediate\write\gtreffile{\\\expandafter\\\def\\\csname ##1\\\endcsname%
{\nextkey}}}}
%
%
\long\def\reflist##1\endreflist{%
\long\def\thereflist{##1}{\def\refkey####1####2\par{\xdef####1{%
\number\refnumber}{\def\\{\noexpand}\immediate\write\gtreffile
{\\\def\\####1{\number\refnumber}}}\global\advance\refnumber by 1}%
\def\key####1####2\par{\expandafter\xdef%
\csname####1\endcsname{\number\refnumber}%
{\def\\{\noexpand}\immediate\write\gtreffile
{\\\expandafter\\\def\\\csname ####1\\\endcsname{\number\refnumber}}}
\global\advance\refnumber by 1}##1\par}}
\long\def\biblio##1\endbiblio{\reflist##1\endreflist\references}%
%
%
\def\numkey##1{{\def\\{\noexpand}%
\xdef##1{\number\sectionnumber.\number\resultnumber}
\immediate\write\gtreffile{\\\def\\##1%
{\number\sectionnumber.\number\resultnumber}}}}
\def\seckey##1{{\def\\{\noexpand}\xdef##1{\number\sectionnumber}
\immediate\write\gtreffile{\\\def\\##1{\number\sectionnumber}}}}
\def\figkey##1{\xdef##1{\number\figurenumber}%
{\def\\{\noexpand}\immediate\write\gtreffile%
{\\\def\\##1{\number\figurenumber}}}
\number\figurenumber\global\advance\figurenumber by 1}
}   
%
%
%
%
\def\figkey#1{\xdef#1{\number\figurenumber}%
\number\figurenumber\global\advance\figurenumber by 1}
\def\fig#1#2\endfig{%
\midinsert\cl{#2}\vglue6pt\cl{\small Figure #1}\endinsert}
\def\newfig{\number\figurenumber\global\advance\figurenumber by 1}
\def\numkey#1{\xdef#1{\number\sectionnumber.\number\resultnumber}}
\def\seckey#1{\xdef#1{\number\sectionnumber}}
%
%
%
%
%
%
%
%
%
\def\verb{\catcode`\"=\active}       
\def\brev{\catcode`\"=12}            
\brev                                
\verb                                
{\obeyspaces\gdef {\ }}              
{\catcode`\`=\active\gdef`{\relax\lq}}
\def"{%
\begingroup\baselineskip=12pt\def\par{\leavevmode\endgraf}%
\tt\obeylines\obeyspaces\parskip=0pt\parindent=0pt%
\catcode`\$=12\catcode`\&=12\catcode`\^=12\catcode`\#=12%
\catcode`\_=12\catcode`\~=12%
\catcode`\{=12\catcode`\}=12\catcode`\%=12\catcode`\\=12%
\catcode`\`=\active\let"\endgroup}
\brev      
%
%
%
%
%
%
\def\item#1{\par\leavevmode\llap{#1\stdspace}%
\ignorespaces}                             
%
%

%
%
\def\co{\colon\thinspace}    
\def\np{\vfil\eject}         
\def\nl{\hfil\break}         
\def\cl{\centerline}         
\def\gt{{\mathsurround=0pt\it $\cal G\mskip-2mu$eometry \&\ 
$\cal T\!\!$opology}}        
\def\agt{{\mathsurround=0pt\it$\cal A\mskip-.7mu$lgebraic \&\ 
$\cal G\mskip-2mu$eometric $\cal T\!\!$opology}}  
%
%
%

%
%
%
%
%
\def\title#1{\def\thetitle{#1}}

\def\author#1{\edef\previousauthors{\theauthors}
 \ifx\theauthors\relax\def\theauthors{#1}\else
 \def\theauthors{\previousauthors\par#1}\fi}

\let\authors\author        
\def\address#1{\edef\previousaddresses{\theaddress}
 \ifx\theaddress\relax\def\theaddress{#1}\else
 \def\theaddress{\previousaddresses\par\vskip 2pt\par#1}\fi}
\def\secondaddress#1{\edef\previousaddresses{\theaddress}
 \ifx\theaddress\relax\def\theaddress{#1}\else
 \def\theaddress{\previousaddresses\par{\rm and}\par#1}\fi}   

\def\email#1{\edef\previousemails{\theemail}
 \ifx\theemail\relax\def\theemail{#1}\else
 \def\theemail{\previousemails\hskip 0.75em\relax#1}\fi}
\def\secondemail#1{\edef\previousemails{\theemail}
 \ifx\theemail\relax\def\theemail{#1}\else
 \def\theemail{\previousemails\hskip 0.75em{\rm and}\hskip 0.75em
 \relax#1}\fi}
\def\url#1{\edef\previousurls{\theurl}
 \ifx\theurl\relax\def\theurl{#1}\else
 \def\theurl{\previousurls\hskip 0.75em\relax#1}\fi}
\def\secondurl#1{\edef\previousurls{\theurl}
 \ifx\theurl\relax\def\theurl{#1}\else
 \def\theurl{\previousurls\hskip 0.75em{\rm and}\hskip 0.75em
 \relax#1}\fi}
\long\def\abstract#1\endabstract{\long\def\theabstract{#1}}
\def\primaryclass#1{\def\theprimaryclass{#1}}
\def\secondaryclass#1{\def\thesecondaryclass{#1}}
\def\keywords#1{\def\thekeywords{#1}}
%
%
\let\\\par\let\thetitle\relax\let\theshorttitle\relax
\let\theauthors\relax\let\theshortauthors\relax
\let\theaddress\relax\let\theshortaddress\relax
\let\theemail\relax\let\theurl\relax
\let\theabstract\relax\let\theprimaryclass\relax
\let\thesecondaryclass\relax\let\thekeywords\relax
%
%
%
%
\long\def\maketitlepage{    

\vglue 0.2truein   

%
{\parskip=0pt\leftskip 0pt plus 1fil\def\\{\par\smallskip}{\large
\bf\thetitle}\par\medskip}   

\vglue 0.15truein 

%
{\parskip=0pt\leftskip 0pt plus 1fil\def\\{\par}{\sc\theauthors}
\par\medskip}%
 
\vglue 0.1truein 

%
{\small\parskip=0pt
{\leftskip 0pt plus 1fil\def\\{\par}{\sl\theaddress}\par}
\ifx\theemail\relax\else  
\vglue 5pt \def\\{\stdspace{\rm and}\stdspace} 
\cl{Email:\stdspace\tt\theemail}\fi
\ifx\theurl\relax\else    
\vglue 5pt \def\\{\stdspace{\rm and}\stdspace} 
\cl{URL:\stdspace\tt\theurl}\fi\par}

\vglue 7pt 

{\bf Abstract}

\vglue 5pt

\theabstract

\vglue 7pt 

{\bf AMS Classification numbers}\quad Primary:\quad \theprimaryclass\par

Secondary:\quad \thesecondaryclass

\vglue 5pt 

{\bf Keywords:}\quad \thekeywords

\np  

}    
%
%
\long\def\makeshorttitle{    


%
{\parskip=0pt\leftskip 0pt plus 1fil\def\\{\par\smallskip}{\large
\bf\thetitle}\par\medskip}   

\vglue 0.05truein 

%
{\parskip=0pt\leftskip 0pt plus 1fil\def\\{\par}{\sc\theauthors}
\par\medskip}%
 
\vglue 0.03truein 

%
{\small\parskip=0pt
{\leftskip 0pt plus 1fil\def\\{\par}{\sl\ifx\theshortaddress\relax
\theaddress\else\theshortaddress\fi}\par}
\ifx\theemail\relax\else  
\vglue 5pt \def\\{\stdspace{\rm and}\stdspace} 
\cl{Email:\stdspace\tt\theemail}\fi
\ifx\theurl\relax\else    
\vglue 5pt \def\\{\stdspace{\rm and}\stdspace} 
\cl{URL:\stdspace\tt\theurl}\fi\par}

\vglue 10pt 


{\small\leftskip 25pt\rightskip 25pt{\bf Abstract}\stdspace\theabstract

{\bf AMS Classification}\stdspace\theprimaryclass
\ifx\thesecondaryclass\relax\else; \thesecondaryclass\fi\par
{\bf Keywords}\stdspace \thekeywords\par}
\vglue 7pt
}    
\let\maketitle\makeshorttitle        
%
%

\def\volumenumber#1{\def\thevolumenumber{#1}}
\def\volumeyear#1{\def\thevolumeyear{#1}}
\def\pagenumbers#1#2{\def\startpage{#1}\def\finishpage{#2}}
\def\published#1{\def\publishdate{#1}}
\def\received#1{\def\receiveddate{#1}}
\def\revised#1{\def\reviseddate{#1}}
\let\reviseddate\relax
\volumenumber{X}
\volumeyear{20XX}
\pagenumbers{1}{XXX}
\published{XX Xxxember 20XX}

\long\def\makeagttitle{   
\agt\hfill      
\hbox to 60truept{\vbox to 0pt{\vglue -14truept{\bf [Logo here]}\vss}\hss}
\break
{\small Volume \thevolumenumber\ (\thevolumeyear)
\startpage--\finishpage\nl
Published: \publishdate}

\vglue .2truein

{\parskip=0pt\leftskip 0pt plus 1fil\def\\{\par\smallskip}{\large
\bf\thetitle}\par\medskip}   
\vglue 0.05truein 

%
{\parskip=0pt\leftskip 0pt plus 1fil\def\\{\par}{\sc\theauthors}
\par\medskip}%
 
\vglue 0.03truein 


{\small\leftskip 25truept\rightskip 25truept{\bf Abstract}\stdspace\theabstract

{\bf AMS Classification}\stdspace\theprimaryclass
\ifx\thesecondaryclass\relax\else; \thesecondaryclass\fi\par
{\bf Keywords}\stdspace \thekeywords\par}\vglue 7truept

}   


\def\Addresses{\bigskip
{\small \parskip 0pt \leftskip 0pt \rightskip 0pt plus 1fil \def\\{\par}
\sl\theaddress\par\medskip \rm Email:\stdspace\tt\theemail\par
\ifx\theurl\relax\else\smallskip \rm URL:\stdspace\tt\theurl\par\fi}}

\def\agtart{
\hoffset 14truemm
\voffset 31truemm
\font\phead=cmsl9 scaled 950
\font\pnum=cmbx10 scaled 913
\font\pfoot=cmsl9 scaled 950
\headline{\vbox to 0pt{\vskip -4.5mm\line{\small\phead\ifnum
\count0=\startpage ISSN numbers are printed here
\hfill {\pnum\folio}\else\ifodd\count0\def\\{ }%
\ifx\theshorttitle\relax\thetitle\else\theshorttitle\fi\hfill{\pnum\folio}
\else\def\\{ and }{\pnum\folio}\hfill\ifx\theshortauthors\relax\theauthors
\else\theshortauthors\fi\fi\fi}\vss}}
\footline{\vbox to 0pt{\vglue 0mm\line{\small\pfoot\ifnum\count0=\startpage
Copyright declaration is printed here\hfill\else
\agt, Volume \thevolumenumber\ (\thevolumeyear)\hfill\fi}\vss}}
\let\maketitle\makeagttitle\let\makeshorttitle\makeagttitle}

\input epsf
%
%
\def\relabelbox{%
  \hbox\bgroup%
}%
\def\endrelabelbox{%
\egroup%
}%
\def\relabel #1#2 {%
  \special{ps:/a {} def}%
  \smash{\rlap{#2}}%
}%
\def\adjustrelabel <#1,#2> #3#4 {%
  \special{ps:/a {} def}%
  \smash{\rlap{\kern #1 \raise #2\hbox{#4}}}%
}%
\def\extralabel <#1,#2> #3 {\smash{\rlap{\kern #1 \raise #2\hbox{#3}}}}%

\chardef\newinsCatAt\the\catcode `\@
\catcode `\@=11
%
%
%
\newskip\insertskipamount\newskip\inserthardskipamount
\insertskipamount 12pt plus2pt     
\inserthardskipamount 4pt          
\def\insertskip{\vskip\insertskipamount}
%
%
\newskip\LastSkip
\def\SaveLastSkip{\LastSkip\lastskip}
\def\RestoreLastSkip{\nobreak\vskip-\LastSkip\vskip\LastSkip}
%
%
\newcount\SplitTest
\def\SetSplitTest{\SplitTest\insertpenalties
  \insert\topins{\floatingpenalty1}%
  \advance\SplitTest-\insertpenalties}
%
%
\def\midinsert{\par
 \SaveLastSkip\penalty-150\SetSplitTest\RestoreLastSkip
 \ifnum\SplitTest=-1
  \@midfalse\p@gefalse\else\@midtrue\fi\@ins}
\def\@ins{\par\begingroup\setbox\z@\vbox\bgroup%
  \vglue\inserthardskipamount}
\def\endinsert{\egroup 
  \if@mid \dimen@\ht\z@ \advance\dimen@\dp\z@
    \advance\dimen@\insertskipamount
    \advance\dimen@\pagetotal\advance\dimen@-\pageshrink
    \ifdim\dimen@>\pagegoal\@midfalse\p@gefalse\fi\fi
  \if@mid%
    \ifdim\lastskip<\insertskipamount\removelastskip\insertskip\fi
    \nointerlineskip\box\z@\penalty-200\insertskip
  \else%
    \SaveLastSkip
    \insert\topins{\penalty100 
    \splittopskip\z@skip
    \splitmaxdepth\maxdimen \floatingpenalty\z@
    \ifp@ge \dimen@\dp\z@
    \vbox to\vsize{\unvbox\z@\kern-\dimen@}
    \else \box\z@\nobreak\insertskip\fi}
    \RestoreLastSkip
   \fi\endgroup}
%
\catcode `\@=\newinsCatAt
 


\def\ifplaintex{\expandafter\ifx\csname documentclass\endcsname\relax}


\ifplaintex 
\hoffset 14truemm
\voffset 31truemm
\else
\headsep 23pt
\footskip 35pt
\hoffset -4truemm
\voffset 12.5truemm
\fi

\expandafter\ifx\csname beginpicture\endcsname\relax
\expandafter\ifx\csname documentclass\endcsname\relax
\input pictex \else\font\fiverm=cmr5
\input prepictex \input pictex \input postpictex \fi\fi

\def\gt{{\mathsurround=0pt\it $\cal G\mskip-2mu$eometry \&\ 
$\cal T\!\!$opology}}        

\def\gtp{{\mathsurround=0pt\it $\cal G\mskip-2mu$eometry \&\ 
$\cal T\!\!$opology $\cal P\!$ublications}}  


\def\lognumber#1{\def\thelognumber{#1}}
\def\volumenumber#1{\def\thevolumenumber{#1}}
\def\papernumber#1{\def\thepapernumber{#1}}
\def\volumeyear#1{\def\thevolumeyear{#1}}

\def\pagenumbers#1#2{\def\startpage{#1}\def\finishpage{#2}}
\def\published#1{\def\publishdate{#1}}
\def\proposed#1{\def\theproposer{#1}}
\def\seconded#1{\def\theseconders{#1}}
\def\received#1{\def\receiveddate{#1}}
\def\revised#1{\def\reviseddate{#1}}
\def\accepted#1{\def\accepteddate{#1}}

\long\def\asciiabstract#1{\long\def\theasciiabstract{#1}}


\let\\\par\let\thelognumber\relax
\let\thevolumenumber\relax\let\thepapernumber\relax
\let\thevolumeyear\relax\let\thesamplenumber\relax\let\startpage\relax
\let\finishpage\relax\let\publishdate\relax\let\receiveddate\relax
\let\reviseddate\relax\let\accepteddate\relax\let\theasciititle\relax
\let\theasciiauthors\relax
\let\theasciiabstract\relax
\let\theasciiemail\relax\let\theshortauthors\relax\let\theshorttitle\relax

\long\def\maketitlep{   

\count0=\startpage

\gt\hfill      
\beginpicture
\setcoordinatesystem units <0.33truein, 0.33truein> point at 2.2 0.9
\setplotsymbol ({$\cal G$})
\plotsymbolspacing=9truept
\circulararc 315 degrees from 0 1 center at 0 0
\setplotsymbol ({$\cal T$})
\circulararc 315 degrees from 1 -1 center at 1 0
\endpicture
%
\break
{\small\ifx\thesamplenumber\relax 
Volume \else Sample
\fi\thevolumenumber\ (\thevolumeyear)
\startpage--\finishpage\nl
Published: \publishdate}
\vglue 0.5truein plus 0.4fil minus 0.1truein

{\parskip=0pt\leftskip 0pt plus 1fil\def\\{\par\smallskip}{\ifplaintex\large
\else\Large\fi\bf\thetitle}\par\medskip}   

\vglue 0pt plus 0.1fil 

{\parskip=0pt\leftskip 0pt plus 1fil\def\\{\par}{\sc\theauthors}
\par\medskip}

\vglue 0pt plus 0.1fil 

{\small\parskip=0pt\let\newline\\
{\leftskip 0pt plus 1fil\def\\{\par}{\sl\theaddress}\par}
\expandafter\ifx\theemail\relax    
\relax\else\vglue 5pt plus 0.02fil minus 2pt\def\\{\stdspace{\rm 
and}\stdspace} 
\cl{Email:\stdspace\tt\theemail}\fi
\ifx\theurl\relax                  
\relax\else\vglue 5pt plus 0.02fil minus 2pt\def\\{\stdspace{\rm 
and}\stdspace}
\cl{URL:\stdspace\tt\theurl}\fi\par}

\vglue 7pt plus 0.3fil minus 3pt

{\bf Abstract}
\vglue 5pt plus 0.1fil minus 2pt

\theabstract

\vglue 7pt plus 0.3fil minus 3pt

{\bf AMS Classification numbers}\quad Primary:\quad \theprimaryclass

Secondary:\quad \thesecondaryclass

\vglue 5pt plus 0.3fil minus 2pt

{\bf Keywords}\quad \thekeywords

\vglue 10pt plus 0.5fil minus 5pt

{\small  Proposed: \theproposer\hfill Received: \receiveddate\nl
Seconded: \theseconders\hfill 
\ifx\reviseddate\relax                         
Accepted: \accepteddate                        
\else
Revised: \reviseddate                          
\fi}
\eject
}       

\let\maketitlepage\maketitlep
\let\maketitle\maketitlepage


\font\phead=cmsl9 scaled 950
\font\lhead=cmsl9 scaled 1050
\font\pnum=cmbx10 scaled 913
\font\lnum=cmbx10 
\font\pfoot=cmsl9 scaled 950
\font\lfoot=cmsl9 scaled 1050
\ifplaintex
\headline{\vbox to 0pt{\vskip -4.5mm\line{\small\phead\ifnum
\count0=\startpage ISSN 1364-0380 (on line)
1465-3060 (printed) \hfill {\pnum\folio}\else\ifodd\count0\def\\{ }%
\ifx\theshorttitle\relax\thetitle\else\theshorttitle\fi\hfill{\pnum\folio}
\else\def\\{ and }{\pnum\folio}\hfill\ifx\theshortauthors\relax\theauthors
\else\theshortauthors\fi\fi\fi}\vss}}
\footline{\vbox to 0pt{\vglue 0mm\line{\small\pfoot\ifnum\count0=\startpage
\copyright\ \gtp\hfill\else
\gt, Volume \thevolumenumber\ (\thevolumeyear)\hfill\fi}\vss
}}
\else
\makeatletter
\def\@oddhead{{\small\lhead\ifnum\count0=\startpage ISSN 1364-0380 (on line)
1465-3060 (printed) \hfill {\lnum\number\count0}\else\ifodd\count0
\def\\{ }\ifx\theshorttitle\relax \thetitle \else\theshorttitle\fi\hfill
{\lnum\number\count0}\else\def\\{ and }{\lnum\number\count0}
\hfill\ifx\theshortauthors\relax 
\theauthors\else\theshortauthors\fi\fi\fi}}\def\@evenhead{\@oddhead}
\def\@oddfoot{\small\lfoot\ifnum\count0=\startpage\copyright\ \gtp\hfill\else
\gt, Volume \thevolumenumber\ (\thevolumeyear)\hfill\fi}
\def\@evenfoot{\@oddfoot}
\makeatother
\fi


\newwrite\gtoutfile
\long\gdef\makeheadfile{  
{\def\\{, }\def\s{ }
\immediate\openout\gtoutfile head.xxx
\immediate\write\gtoutfile{To: math@arxiv.org}
\immediate\write\gtoutfile{Subject: put or rep NNNNN:pppp}
\immediate\write\gtoutfile{--text follows this line--}
\immediate\write\gtoutfile{Proxy-for: \ifx\theasciiauthors\relax
\theauthors\else\theasciiauthors\fi\s<\ifx\theasciiemail\relax\theemail\else\theasciiemail\fi>}
\immediate\write\gtoutfile{\noexpand\\}
\immediate\write\gtoutfile{Authors: \ifx\theasciiauthors\relax
\theauthors\else\theasciiauthors\fi}
\immediate\write\gtoutfile{Title: \ifx\theasciititle\relax
\thetitle\else\theasciititle\fi}
\immediate\write\gtoutfile{Subj-class: GT or SG or MG etc}
\immediate\write\gtoutfile{MSC-class: \theprimaryclass\ifx\thesecondaryclass\relax\else, \thesecondaryclass\fi}
\immediate\write\gtoutfile{Journal-ref: Geom. Topol. \thevolumenumber
(\thevolumeyear) \startpage-\finishpage}
\immediate\write\gtoutfile{Comments: Published by Geometry and Topology at}
\immediate\write\gtoutfile{\s\s http://www.maths.warwick.ac.uk/gt/GTVol\thevolumenumber/paper\thepapernumber.abs.html}
\immediate\write\gtoutfile{\noexpand\\}
\immediate\write\gtoutfile{}
\ifx\theasciiabstract\relax
\immediate\write\gtoutfile{\theabstract}\else
\immediate\write\gtoutfile{\theasciiabstract}\fi
\immediate\write\gtoutfile{}
\immediate\write\gtoutfile{\noexpand\\}
\immediate\write\gtoutfile{}
\immediate\closeout\gtoutfile}}  

\def\maketitlepage{\maketitlep\makeheadfile}
\let\maketitle\maketitlepage


\def\ifplaintex{\expandafter\ifx\csname documentclass\endcsname\relax}


\ifplaintex 
\hoffset 14truemm
\voffset 31truemm
\else
\headsep 23pt
\footskip 35pt
\hoffset -4truemm
\voffset 12.5truemm
\fi

\expandafter\ifx\csname beginpicture\endcsname\relax
\expandafter\ifx\csname documentclass\endcsname\relax
\input pictex \else\font\fiverm=cmr5
\input prepictex \input pictex \input postpictex \fi\fi

\def\gt{{\mathsurround=0pt\it $\cal G\mskip-2mu$eometry \&\ 
$\cal T\!\!$opology}}        

\def\gtp{{\mathsurround=0pt\it $\cal G\mskip-2mu$eometry \&\ 
$\cal T\!\!$opology $\cal P\!$ublications}}  


\def\lognumber#1{\def\thelognumber{#1}}
\def\volumenumber#1{\def\thevolumenumber{#1}}
\def\papernumber#1{\def\thepapernumber{#1}}
\def\volumeyear#1{\def\thevolumeyear{#1}}

\def\pagenumbers#1#2{\def\startpage{#1}\def\finishpage{#2}}
\def\published#1{\def\publishdate{#1}}
\def\proposed#1{\def\theproposer{#1}}
\def\seconded#1{\def\theseconders{#1}}
\def\received#1{\def\receiveddate{#1}}
\def\revised#1{\def\reviseddate{#1}}
\def\accepted#1{\def\accepteddate{#1}}

\long\def\asciiabstract#1{\long\def\theasciiabstract{#1}}


\let\\\par\let\thelognumber\relax
\let\thevolumenumber\relax\let\thepapernumber\relax
\let\thevolumeyear\relax\let\thesamplenumber\relax\let\startpage\relax
\let\finishpage\relax\let\publishdate\relax\let\receiveddate\relax
\let\reviseddate\relax\let\accepteddate\relax\let\theasciititle\relax
\let\theasciiauthors\relax
\let\theasciiabstract\relax
\let\theasciiemail\relax\let\theshortauthors\relax\let\theshorttitle\relax

\long\def\maketitlep{   

\count0=\startpage

\gt\hfill      
\beginpicture
\setcoordinatesystem units <0.33truein, 0.33truein> point at 2.2 0.9
\setplotsymbol ({$\cal G$})
\plotsymbolspacing=9truept
\circulararc 315 degrees from 0 1 center at 0 0
\setplotsymbol ({$\cal T$})
\circulararc 315 degrees from 1 -1 center at 1 0
\endpicture
%
\break
{\small\ifx\thesamplenumber\relax 
Volume \else Sample
\fi\thevolumenumber\ (\thevolumeyear)
\startpage--\finishpage\nl
Published: \publishdate}
\vglue 0.5truein plus 0.4fil minus 0.1truein

{\parskip=0pt\leftskip 0pt plus 1fil\def\\{\par\smallskip}{\ifplaintex\large
\else\Large\fi\bf\thetitle}\par\medskip}   

\vglue 0pt plus 0.1fil 

{\parskip=0pt\leftskip 0pt plus 1fil\def\\{\par}{\sc\theauthors}
\par\medskip}

\vglue 0pt plus 0.1fil 

{\small\parskip=0pt\let\newline\\
{\leftskip 0pt plus 1fil\def\\{\par}{\sl\theaddress}\par}
\expandafter\ifx\theemail\relax    
\relax\else\vglue 5pt plus 0.02fil minus 2pt\def\\{\stdspace{\rm 
and}\stdspace} 
\cl{Email:\stdspace\tt\theemail}\fi
\ifx\theurl\relax                  
\relax\else\vglue 5pt plus 0.02fil minus 2pt\def\\{\stdspace{\rm 
and}\stdspace}
\cl{URL:\stdspace\tt\theurl}\fi\par}

\vglue 7pt plus 0.3fil minus 3pt

{\bf Abstract}
\vglue 5pt plus 0.1fil minus 2pt

\theabstract

\vglue 7pt plus 0.3fil minus 3pt

{\bf AMS Classification numbers}\quad Primary:\quad \theprimaryclass

Secondary:\quad \thesecondaryclass

\vglue 5pt plus 0.3fil minus 2pt

{\bf Keywords}\quad \thekeywords

\vglue 10pt plus 0.5fil minus 5pt

{\small  Proposed: \theproposer\hfill Received: \receiveddate\nl
Seconded: \theseconders\hfill 
\ifx\reviseddate\relax                         
Accepted: \accepteddate                        
\else
Revised: \reviseddate                          
\fi}
\eject
}       

\let\maketitlepage\maketitlep
\let\maketitle\maketitlepage


\font\phead=cmsl9 scaled 950
\font\lhead=cmsl9 scaled 1050
\font\pnum=cmbx10 scaled 913
\font\lnum=cmbx10 
\font\pfoot=cmsl9 scaled 950
\font\lfoot=cmsl9 scaled 1050
\ifplaintex
\headline{\vbox to 0pt{\vskip -4.5mm\line{\small\phead\ifnum
\count0=\startpage ISSN 1364-0380 (on line)
1465-3060 (printed) \hfill {\pnum\folio}\else\ifodd\count0\def\\{ }%
\ifx\theshorttitle\relax\thetitle\else\theshorttitle\fi\hfill{\pnum\folio}
\else\def\\{ and }{\pnum\folio}\hfill\ifx\theshortauthors\relax\theauthors
\else\theshortauthors\fi\fi\fi}\vss}}
\footline{\vbox to 0pt{\vglue 0mm\line{\small\pfoot\ifnum\count0=\startpage
\copyright\ \gtp\hfill\else
\gt, Volume \thevolumenumber\ (\thevolumeyear)\hfill\fi}\vss
}}
\else
\makeatletter
\def\@oddhead{{\small\lhead\ifnum\count0=\startpage ISSN 1364-0380 (on line)
1465-3060 (printed) \hfill {\lnum\number\count0}\else\ifodd\count0
\def\\{ }\ifx\theshorttitle\relax \thetitle \else\theshorttitle\fi\hfill
{\lnum\number\count0}\else\def\\{ and }{\lnum\number\count0}
\hfill\ifx\theshortauthors\relax 
\theauthors\else\theshortauthors\fi\fi\fi}}\def\@evenhead{\@oddhead}
\def\@oddfoot{\small\lfoot\ifnum\count0=\startpage\copyright\ \gtp\hfill\else
\gt, Volume \thevolumenumber\ (\thevolumeyear)\hfill\fi}
\def\@evenfoot{\@oddfoot}
\makeatother
\fi


\newwrite\gtoutfile
\long\gdef\makeheadfile{  
{\def\\{, }\def\s{ }
\immediate\openout\gtoutfile head.xxx
\immediate\write\gtoutfile{To: math@arxiv.org}
\immediate\write\gtoutfile{Subject: put or rep NNNNN:pppp}
\immediate\write\gtoutfile{--text follows this line--}
\immediate\write\gtoutfile{Proxy-for: \ifx\theasciiauthors\relax
\theauthors\else\theasciiauthors\fi\s<\ifx\theasciiemail\relax\theemail\else\theasciiemail\fi>}
\immediate\write\gtoutfile{\noexpand\\}
\immediate\write\gtoutfile{Authors: \ifx\theasciiauthors\relax
\theauthors\else\theasciiauthors\fi}
\immediate\write\gtoutfile{Title: \ifx\theasciititle\relax
\thetitle\else\theasciititle\fi}
\immediate\write\gtoutfile{Subj-class: GT or SG or MG etc}
\immediate\write\gtoutfile{MSC-class: \theprimaryclass\ifx\thesecondaryclass\relax\else, \thesecondaryclass\fi}
\immediate\write\gtoutfile{Journal-ref: Geom. Topol. \thevolumenumber
(\thevolumeyear) \startpage-\finishpage}
\immediate\write\gtoutfile{Comments: Published by Geometry and Topology at}
\immediate\write\gtoutfile{\s\s http://www.maths.warwick.ac.uk/gt/GTVol\thevolumenumber/paper\thepapernumber.abs.html}
\immediate\write\gtoutfile{\noexpand\\}
\immediate\write\gtoutfile{}
\ifx\theasciiabstract\relax
\immediate\write\gtoutfile{\theabstract}\else
\immediate\write\gtoutfile{\theasciiabstract}\fi
\immediate\write\gtoutfile{}
\immediate\write\gtoutfile{\noexpand\\}
\immediate\write\gtoutfile{}
\immediate\closeout\gtoutfile}}  

\def\maketitlepage{\maketitlep\makeheadfile}
\let\maketitle\maketitlepage


\def\ifplaintex{\expandafter\ifx\csname documentclass\endcsname\relax}


\ifplaintex 
\hoffset 14truemm
\voffset 31truemm
\else
\headsep 23pt
\footskip 35pt
\hoffset -4truemm
\voffset 12.5truemm
\fi

\expandafter\ifx\csname beginpicture\endcsname\relax
\expandafter\ifx\csname documentclass\endcsname\relax
\input pictex \else\font\fiverm=cmr5
\input prepictex \input pictex \input postpictex \fi\fi

\def\gt{{\mathsurround=0pt\it $\cal G\mskip-2mu$eometry \&\ 
$\cal T\!\!$opology}}        

\def\gtp{{\mathsurround=0pt\it $\cal G\mskip-2mu$eometry \&\ 
$\cal T\!\!$opology $\cal P\!$ublications}}  


\def\lognumber#1{\def\thelognumber{#1}}
\def\volumenumber#1{\def\thevolumenumber{#1}}
\def\papernumber#1{\def\thepapernumber{#1}}
\def\volumeyear#1{\def\thevolumeyear{#1}}

\def\pagenumbers#1#2{\def\startpage{#1}\def\finishpage{#2}}
\def\published#1{\def\publishdate{#1}}
\def\proposed#1{\def\theproposer{#1}}
\def\seconded#1{\def\theseconders{#1}}
\def\received#1{\def\receiveddate{#1}}
\def\revised#1{\def\reviseddate{#1}}
\def\accepted#1{\def\accepteddate{#1}}

\long\def\asciiabstract#1{\long\def\theasciiabstract{#1}}


\let\\\par\let\thelognumber\relax
\let\thevolumenumber\relax\let\thepapernumber\relax
\let\thevolumeyear\relax\let\thesamplenumber\relax\let\startpage\relax
\let\finishpage\relax\let\publishdate\relax\let\receiveddate\relax
\let\reviseddate\relax\let\accepteddate\relax\let\theasciititle\relax
\let\theasciiauthors\relax
\let\theasciiabstract\relax
\let\theasciiemail\relax\let\theshortauthors\relax\let\theshorttitle\relax

\long\def\maketitlep{   

\count0=\startpage

\gt\hfill      
\beginpicture
\setcoordinatesystem units <0.33truein, 0.33truein> point at 2.2 0.9
\setplotsymbol ({$\cal G$})
\plotsymbolspacing=9truept
\circulararc 315 degrees from 0 1 center at 0 0
\setplotsymbol ({$\cal T$})
\circulararc 315 degrees from 1 -1 center at 1 0
\endpicture
%
\break
{\small\ifx\thesamplenumber\relax 
Volume \else Sample
\fi\thevolumenumber\ (\thevolumeyear)
\startpage--\finishpage\nl
Published: \publishdate}
\vglue 0.5truein plus 0.4fil minus 0.1truein

{\parskip=0pt\leftskip 0pt plus 1fil\def\\{\par\smallskip}{\ifplaintex\large
\else\Large\fi\bf\thetitle}\par\medskip}   

\vglue 0pt plus 0.1fil 

{\parskip=0pt\leftskip 0pt plus 1fil\def\\{\par}{\sc\theauthors}
\par\medskip}

\vglue 0pt plus 0.1fil 

{\small\parskip=0pt\let\newline\\
{\leftskip 0pt plus 1fil\def\\{\par}{\sl\theaddress}\par}
\expandafter\ifx\theemail\relax    
\relax\else\vglue 5pt plus 0.02fil minus 2pt\def\\{\stdspace{\rm 
and}\stdspace} 
\cl{Email:\stdspace\tt\theemail}\fi
\ifx\theurl\relax                  
\relax\else\vglue 5pt plus 0.02fil minus 2pt\def\\{\stdspace{\rm 
and}\stdspace}
\cl{URL:\stdspace\tt\theurl}\fi\par}

\vglue 7pt plus 0.3fil minus 3pt

{\bf Abstract}
\vglue 5pt plus 0.1fil minus 2pt

\theabstract

\vglue 7pt plus 0.3fil minus 3pt

{\bf AMS Classification numbers}\quad Primary:\quad \theprimaryclass

Secondary:\quad \thesecondaryclass

\vglue 5pt plus 0.3fil minus 2pt

{\bf Keywords}\quad \thekeywords

\vglue 10pt plus 0.5fil minus 5pt

{\small  Proposed: \theproposer\hfill Received: \receiveddate\nl
Seconded: \theseconders\hfill 
\ifx\reviseddate\relax                         
Accepted: \accepteddate                        
\else
Revised: \reviseddate                          
\fi}
\eject
}       

\let\maketitlepage\maketitlep
\let\maketitle\maketitlepage


\font\phead=cmsl9 scaled 950
\font\lhead=cmsl9 scaled 1050
\font\pnum=cmbx10 scaled 913
\font\lnum=cmbx10 
\font\pfoot=cmsl9 scaled 950
\font\lfoot=cmsl9 scaled 1050
\ifplaintex
\headline{\vbox to 0pt{\vskip -4.5mm\line{\small\phead\ifnum
\count0=\startpage ISSN 1364-0380 (on line)
1465-3060 (printed) \hfill {\pnum\folio}\else\ifodd\count0\def\\{ }%
\ifx\theshorttitle\relax\thetitle\else\theshorttitle\fi\hfill{\pnum\folio}
\else\def\\{ and }{\pnum\folio}\hfill\ifx\theshortauthors\relax\theauthors
\else\theshortauthors\fi\fi\fi}\vss}}
\footline{\vbox to 0pt{\vglue 0mm\line{\small\pfoot\ifnum\count0=\startpage
\copyright\ \gtp\hfill\else
\gt, Volume \thevolumenumber\ (\thevolumeyear)\hfill\fi}\vss
}}
\else
\makeatletter
\def\@oddhead{{\small\lhead\ifnum\count0=\startpage ISSN 1364-0380 (on line)
1465-3060 (printed) \hfill {\lnum\number\count0}\else\ifodd\count0
\def\\{ }\ifx\theshorttitle\relax \thetitle \else\theshorttitle\fi\hfill
{\lnum\number\count0}\else\def\\{ and }{\lnum\number\count0}
\hfill\ifx\theshortauthors\relax 
\theauthors\else\theshortauthors\fi\fi\fi}}\def\@evenhead{\@oddhead}
\def\@oddfoot{\small\lfoot\ifnum\count0=\startpage\copyright\ \gtp\hfill\else
\gt, Volume \thevolumenumber\ (\thevolumeyear)\hfill\fi}
\def\@evenfoot{\@oddfoot}
\makeatother
\fi


\newwrite\gtoutfile
\long\gdef\makeheadfile{  
{\def\\{, }\def\s{ }
\immediate\openout\gtoutfile head.xxx
\immediate\write\gtoutfile{To: math@arxiv.org}
\immediate\write\gtoutfile{Subject: put or rep NNNNN:pppp}
\immediate\write\gtoutfile{--text follows this line--}
\immediate\write\gtoutfile{Proxy-for: \ifx\theasciiauthors\relax
\theauthors\else\theasciiauthors\fi\s<\ifx\theasciiemail\relax\theemail\else\theasciiemail\fi>}
\immediate\write\gtoutfile{\noexpand\\}
\immediate\write\gtoutfile{Authors: \ifx\theasciiauthors\relax
\theauthors\else\theasciiauthors\fi}
\immediate\write\gtoutfile{Title: \ifx\theasciititle\relax
\thetitle\else\theasciititle\fi}
\immediate\write\gtoutfile{Subj-class: GT or SG or MG etc}
\immediate\write\gtoutfile{MSC-class: \theprimaryclass\ifx\thesecondaryclass\relax\else, \thesecondaryclass\fi}
\immediate\write\gtoutfile{Journal-ref: Geom. Topol. \thevolumenumber
(\thevolumeyear) \startpage-\finishpage}
\immediate\write\gtoutfile{Comments: Published by Geometry and Topology at}
\immediate\write\gtoutfile{\s\s http://www.maths.warwick.ac.uk/gt/GTVol\thevolumenumber/paper\thepapernumber.abs.html}
\immediate\write\gtoutfile{\noexpand\\}
\immediate\write\gtoutfile{}
\ifx\theasciiabstract\relax
\immediate\write\gtoutfile{\theabstract}\else
\immediate\write\gtoutfile{\theasciiabstract}\fi
\immediate\write\gtoutfile{}
\immediate\write\gtoutfile{\noexpand\\}
\immediate\write\gtoutfile{}
\immediate\closeout\gtoutfile}}  

\def\maketitlepage{\maketitlep\makeheadfile}
\let\maketitle\maketitlepage


\def\ifplaintex{\expandafter\ifx\csname documentclass\endcsname\relax}


\ifplaintex 
\hoffset 14truemm
\voffset 31truemm
\else
\headsep 23pt
\footskip 35pt
\hoffset -4truemm
\voffset 12.5truemm
\fi

\expandafter\ifx\csname beginpicture\endcsname\relax
\expandafter\ifx\csname documentclass\endcsname\relax
\input pictex \else\font\fiverm=cmr5
\input prepictex \input pictex \input postpictex \fi\fi

\def\gt{{\mathsurround=0pt\it $\cal G\mskip-2mu$eometry \&\ 
$\cal T\!\!$opology}}        

\def\gtp{{\mathsurround=0pt\it $\cal G\mskip-2mu$eometry \&\ 
$\cal T\!\!$opology $\cal P\!$ublications}}  


\def\lognumber#1{\def\thelognumber{#1}}
\def\volumenumber#1{\def\thevolumenumber{#1}}
\def\papernumber#1{\def\thepapernumber{#1}}
\def\volumeyear#1{\def\thevolumeyear{#1}}

\def\pagenumbers#1#2{\def\startpage{#1}\def\finishpage{#2}}
\def\published#1{\def\publishdate{#1}}
\def\proposed#1{\def\theproposer{#1}}
\def\seconded#1{\def\theseconders{#1}}
\def\received#1{\def\receiveddate{#1}}
\def\revised#1{\def\reviseddate{#1}}
\def\accepted#1{\def\accepteddate{#1}}

\long\def\asciiabstract#1{\long\def\theasciiabstract{#1}}


\let\\\par\let\thelognumber\relax
\let\thevolumenumber\relax\let\thepapernumber\relax
\let\thevolumeyear\relax\let\thesamplenumber\relax\let\startpage\relax
\let\finishpage\relax\let\publishdate\relax\let\receiveddate\relax
\let\reviseddate\relax\let\accepteddate\relax\let\theasciititle\relax
\let\theasciiauthors\relax
\let\theasciiabstract\relax
\let\theasciiemail\relax\let\theshortauthors\relax\let\theshorttitle\relax

\long\def\maketitlep{   

\count0=\startpage

\gt\hfill      
\beginpicture
\setcoordinatesystem units <0.33truein, 0.33truein> point at 2.2 0.9
\setplotsymbol ({$\cal G$})
\plotsymbolspacing=9truept
\circulararc 315 degrees from 0 1 center at 0 0
\setplotsymbol ({$\cal T$})
\circulararc 315 degrees from 1 -1 center at 1 0
\endpicture
%
\break
{\small\ifx\thesamplenumber\relax 
Volume \else Sample
\fi\thevolumenumber\ (\thevolumeyear)
\startpage--\finishpage\nl
Published: \publishdate}
\vglue 0.5truein plus 0.4fil minus 0.1truein

{\parskip=0pt\leftskip 0pt plus 1fil\def\\{\par\smallskip}{\ifplaintex\large
\else\Large\fi\bf\thetitle}\par\medskip}   

\vglue 0pt plus 0.1fil 

{\parskip=0pt\leftskip 0pt plus 1fil\def\\{\par}{\sc\theauthors}
\par\medskip}

\vglue 0pt plus 0.1fil 

{\small\parskip=0pt\let\newline\\
{\leftskip 0pt plus 1fil\def\\{\par}{\sl\theaddress}\par}
\expandafter\ifx\theemail\relax    
\relax\else\vglue 5pt plus 0.02fil minus 2pt\def\\{\stdspace{\rm 
and}\stdspace} 
\cl{Email:\stdspace\tt\theemail}\fi
\ifx\theurl\relax                  
\relax\else\vglue 5pt plus 0.02fil minus 2pt\def\\{\stdspace{\rm 
and}\stdspace}
\cl{URL:\stdspace\tt\theurl}\fi\par}

\vglue 7pt plus 0.3fil minus 3pt

{\bf Abstract}
\vglue 5pt plus 0.1fil minus 2pt

\theabstract

\vglue 7pt plus 0.3fil minus 3pt

{\bf AMS Classification numbers}\quad Primary:\quad \theprimaryclass

Secondary:\quad \thesecondaryclass

\vglue 5pt plus 0.3fil minus 2pt

{\bf Keywords}\quad \thekeywords

\vglue 10pt plus 0.5fil minus 5pt

{\small  Proposed: \theproposer\hfill Received: \receiveddate\nl
Seconded: \theseconders\hfill 
\ifx\reviseddate\relax                         
Accepted: \accepteddate                        
\else
Revised: \reviseddate                          
\fi}
\eject
}       

\let\maketitlepage\maketitlep
\let\maketitle\maketitlepage


\font\phead=cmsl9 scaled 950
\font\lhead=cmsl9 scaled 1050
\font\pnum=cmbx10 scaled 913
\font\lnum=cmbx10 
\font\pfoot=cmsl9 scaled 950
\font\lfoot=cmsl9 scaled 1050
\ifplaintex
\headline{\vbox to 0pt{\vskip -4.5mm\line{\small\phead\ifnum
\count0=\startpage ISSN 1364-0380 (on line)
1465-3060 (printed) \hfill {\pnum\folio}\else\ifodd\count0\def\\{ }%
\ifx\theshorttitle\relax\thetitle\else\theshorttitle\fi\hfill{\pnum\folio}
\else\def\\{ and }{\pnum\folio}\hfill\ifx\theshortauthors\relax\theauthors
\else\theshortauthors\fi\fi\fi}\vss}}
\footline{\vbox to 0pt{\vglue 0mm\line{\small\pfoot\ifnum\count0=\startpage
\copyright\ \gtp\hfill\else
\gt, Volume \thevolumenumber\ (\thevolumeyear)\hfill\fi}\vss
}}
\else
\makeatletter
\def\@oddhead{{\small\lhead\ifnum\count0=\startpage ISSN 1364-0380 (on line)
1465-3060 (printed) \hfill {\lnum\number\count0}\else\ifodd\count0
\def\\{ }\ifx\theshorttitle\relax \thetitle \else\theshorttitle\fi\hfill
{\lnum\number\count0}\else\def\\{ and }{\lnum\number\count0}
\hfill\ifx\theshortauthors\relax 
\theauthors\else\theshortauthors\fi\fi\fi}}\def\@evenhead{\@oddhead}
\def\@oddfoot{\small\lfoot\ifnum\count0=\startpage\copyright\ \gtp\hfill\else
\gt, Volume \thevolumenumber\ (\thevolumeyear)\hfill\fi}
\def\@evenfoot{\@oddfoot}
\makeatother
\fi


\newwrite\gtoutfile
\long\gdef\makeheadfile{  
{\def\\{, }\def\s{ }
\immediate\openout\gtoutfile head.xxx
\immediate\write\gtoutfile{To: math@arxiv.org}
\immediate\write\gtoutfile{Subject: put or rep NNNNN:pppp}
\immediate\write\gtoutfile{--text follows this line--}
\immediate\write\gtoutfile{Proxy-for: \ifx\theasciiauthors\relax
\theauthors\else\theasciiauthors\fi\s<\ifx\theasciiemail\relax\theemail\else\theasciiemail\fi>}
\immediate\write\gtoutfile{\noexpand\\}
\immediate\write\gtoutfile{Authors: \ifx\theasciiauthors\relax
\theauthors\else\theasciiauthors\fi}
\immediate\write\gtoutfile{Title: \ifx\theasciititle\relax
\thetitle\else\theasciititle\fi}
\immediate\write\gtoutfile{Subj-class: GT or SG or MG etc}
\immediate\write\gtoutfile{MSC-class: \theprimaryclass\ifx\thesecondaryclass\relax\else, \thesecondaryclass\fi}
\immediate\write\gtoutfile{Journal-ref: Geom. Topol. \thevolumenumber
(\thevolumeyear) \startpage-\finishpage}
\immediate\write\gtoutfile{Comments: Published by Geometry and Topology at}
\immediate\write\gtoutfile{\s\s http://www.maths.warwick.ac.uk/gt/GTVol\thevolumenumber/paper\thepapernumber.abs.html}
\immediate\write\gtoutfile{\noexpand\\}
\immediate\write\gtoutfile{}
\ifx\theasciiabstract\relax
\immediate\write\gtoutfile{\theabstract}\else
\immediate\write\gtoutfile{\theasciiabstract}\fi
\immediate\write\gtoutfile{}
\immediate\write\gtoutfile{\noexpand\\}
\immediate\write\gtoutfile{}
\immediate\closeout\gtoutfile}}  

\def\maketitlepage{\maketitlep\makeheadfile}
\let\maketitle\maketitlepage

\lognumber{179}
\volumenumber{5}\papernumber{15}\volumeyear{2001}
\pagenumbers{431}{440}
\accepted{23 April 2001}
\received{25 January 2001}
\revised{2 April 2001}
\published{24 April 2001}
\proposed{Robion Kirby}
\seconded{Yasha Eliashberg, David Gabai}
%
%
\def\tw{{\tt \char'176}}  

\def\ep{\varepsilon}

\font\spec=cmtex10  
\def\d{\hbox{\spec \char'017\kern 0.05em}} 
\def\inv{^{-1}}

\reflist

\key{Grom}
{\bf M Gromov}, {\it Partial differential relations}, Springer--Verlag
(1986)

\key{Kuiper}
{\bf Nicolaas H Kuiper}, {\it On $C^1$--isometric imbeddings. I, II},
Nederl. Akad. Wetensch. Proc. Ser. A. {58} = Indag. Math. 17 (1955),
545--556 and 683--689

\key{CompI}
{\bf Colin Rourke}, {\bf Brian Sanderson}, {\it The compression
theorem I}, Geometry and Topology, 5 (2001) 399--429,
{\tt arxiv:math.GT/9712235}

\refkey\CompIII
{\bf Colin Rourke}, {\bf Brian Sanderson}, {\it The compression theorem III:
applications}, to appear, see {\tt http://www.maths.warwick.ac.uk/\tw
cpr/ftp/compIII.ps} for a preliminary version

\endreflist

\title{The compression theorem II: directed embeddings}

\authors{Colin Rourke\\Brian Sanderson}

\address{Mathematics Institute, University of Warwick\\Coventry, CV4 7AL, UK}

\email{cpr@maths.warwick.ac.uk\\bjs@maths.warwick.ac.uk}
\url{http://www.maths.warwick.ac.uk/\tw cpr/\qua {\rm and}\qua \tw bjs/}

\abstract

This is the second of three papers about the Compression Theorem.  We
give proofs of Gromov's theorem on directed embeddings [\Grom; 2.4.5
$\rm(C')$] and of the Normal Deformation Theorem [\CompI; 4.7] (a
general version of the Compression Theorem).

\endabstract

\asciiabstract{
This is the second of three papers about the Compression Theorem.  We
give proofs of Gromov's theorem on directed embeddings [M Gromov,
Partial differential relations, Springer--Verlag (1986); 2.4.5 C'] and
of the Normal Deformation Theorem [The compression theorem I; 4.7] (a
general version of the Compression Theorem).}

\primaryclass{57R40, 57R42}
\secondaryclass{57A05}

\keywords{Compression, flattening, directed embedding, ripple lemma}

\maketitlepage

\section{Introduction}

This is the second of three papers about the Compression Theorem.  The
first paper [\CompI] contains a proof of the theorem using an explicit
vector field argument.  This paper contains one simple piece of
geometry, {\it rippling}, which proves the ``flattening lemma'' stated
below.  This leads to proofs of Gromov's theorem on directed
embeddings [\Grom; 2.4.5 $\rm(C')]$ and of the Normal Deformation
Theorem [\CompI; 4.7] (a general version of the Compression Theorem
which can be readily deduced from Gromov's theorem).  The third paper
[\CompIII] is concerned with applications.

We work throughout in the smooth ($C^\infty$) category and we shall
assume without comment that all manifolds are equipped with
appropriate Riemannian metrics.  The tangent bundle of a manifold $W$
is denoted $TW$.  Suppose that $M^n$ is smoothly embedded in
$Q^q\times\re^t$ where $q\ge n$.  We think of $\re^t$ as {\it
vertical} and $Q$ as {\it horizontal}.  We say that $M$ is {\sl
compressible} if it is nowhere tangent to vertical, or equivalently,
if projection $p$ on $Q$ takes $M$ to an immersion in $Q$.  Throughout
the paper, ``normal'' means independent (as in the usual meaning of
``normal bundle'') and not necessarily perpendicular.

Let $G=G_n(Q^q\times\re^t)$ denote the Grassmann bundle of
$n$--planes in $T(Q^q\times\re^t)$ and define the {\sl horizontal}
subset $H$ of $G$ to comprise $n$--planes with no vertical component.
In other words $H$ comprises $n$--planes lying in fibres of $p^*TQ$.

\proclaim{Flattening lemma} Suppose that $M^n$ is compressible in 
$Q^q\times\re^t$ and that $U$ is any neighbourhood of $H$ in $G$.  Suppose
that $q-n\ge 1$ then there is a $C^0$--small isotopy of $M$ in $Q$
carrying $M$ to a position where $TM\subset U$.\endproc

We think of planes in $H$ as {\sl flat} and planes in $U$ as {\sl
almost} flat.  So the lemma moves $M$ to a position where it is almost
flat (ie, its tangent bundle comprises almost flat planes).  Obviously
it is in general impossible to move $M$ to a position where it is
completely flat.

\rk{Addenda}

(1)\qua The lemma is also true if $q=n$ and each component of $M$ is
either open or has boundary.  However in this case the isotopy is not
$C^0$--small, but of the form ``shrink to a neighbourhood of a
(chosen) spine of $M$ followed by a small isotopy''.  We call such an
isotopy ``pseudo-small''; a pseudo-small isotopy has support in a
small neighbourhood of the original embedding.

(2)\qua The lemma has both relative and parametrised versions: If $M$
is already almost flat in the neighbourhood of some closed subset $C$
of $M$ then the isotopy can be assumed fixed on $C$.  Further given a
family of compressible embeddings of $M$ in $Q$ parametrised by a
manifold $K$, then there is a $K$--parametrised family of small
isotopies of $M$ in $Q$ carrying $TM\times K$ into a given
neighbourhood of $H\times K$.

(3)\qua In the non-compact case smallness can be assumed to vary.  In
other words, during the isotopy, points move no further than $\ep>0$,
which is a given function of $M$ (or $M\times K$ in the parametrised
version).

\section{The proof of the flattening lemma}

The process is analogous to the way in which a road is constructed to
go up a steep hill.  Hairpin bends are added (which has the effect of
greatly increasing the horizontal distance the road travels) and this
allows the slope of the road to become as small as necessary.

The process of adding hairpin bends is embodied in the ``ripple
lemma'' which we state and prove first.  The flattening lemma follows
quickly from the ripple lemma.

Let $M$ be connected and smoothly immersed in $Q$.  Suppose that $Q$
has a Riemannian metric $d$.  Define the {\sl induced Riemannian
metric} denoted $d_M(x,y)$ on $M$ by restricting to $TM$ the form on
$TQ$ which defines $d$.

If $M$ is embedded in $Q$ then we also have the usual induced metric
on $M$ (ie $d$ restricted to $M$).  The two induced metrics coincide
to first order for nearby points but in general the induced metric is
{\it smaller} than the induced Riemannian metric.

\proclaim{Ripple lemma}Suppose that $M$ is smoothly embedded in $Q$
and that $q-n\ge 1$ and that $R,\ep>0$ are any given real
numbers.  Then there is an isotopy of $M$ in $Q$ which moves
points at most $\ep$ such that the finishing diffeomorphism
$f\co M\to f(M)$ has the following property:\smallskip

\line{\hfill$d_{f(M)}(f(x),f(y))>Rd_M(x,y)$ for all $x,y\in M$
\hfill$(*)$}\endproc

Thus the ripple lemma asserts that we can (by a small isotopy) arrange
for distances of points (measured using the induced Riemannian metric
from $Q$) to be scaled up by as large a factor as we please.  The
proof is to systematically ``ripple'' the embedding, hence the name.

\prf 
We deal first with the case when $Q$ is the plane and $M$ is the unit
interval $[0,1]$ in the $x$--axis, so $n=1$ and $q=2$.  We will work
relative to the boundary of $M$; in other words the isotopy we
construct will be fixed near $0$ and $1$ and the scaling up will work
for points outside a given neighbourhood of $\{0,1\}$.

Consider a sine curve $S$ of amplitude $A$ and frequency $\omega$ (the
graph of $y=A\sin(2\pi\omega x)$).  Think of $\omega$ as large and $A$
as small.  So the curve is a small ripple of high frequency.  Use a
$C^\infty$ bump function to phase $S$ down to zero near $x=0$ and
$x=1$.  See figure 1.

\figure\relabelbox\small
\epsfxsize 3 truein
\epsfbox{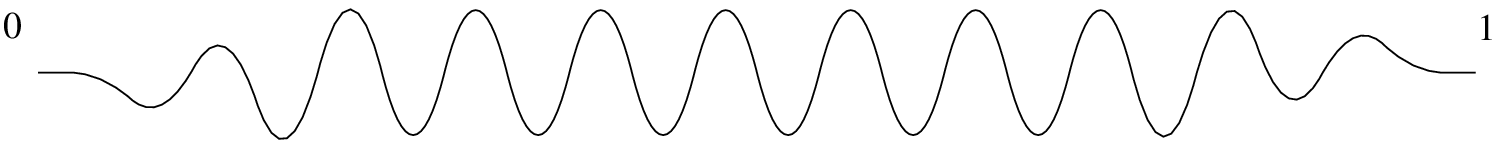}
\relabel {0}{$0$}
\relabel {1}{$1$}
\endrelabelbox
\caption{The basic ripple}
\endfigure

Clearly $M$ can be replaced by $S$ via an isotopy which moves points
at most $A+{1\over \omega}$.  Further we can choose this isotopy to
finish with a diffeomorphism which nowhere shrinks distances, is fixed
near $\{0,1\}$, and which outside a given neighbourhood of $\{0,1\}$
scales length up (measured in $M$) by a constant scale factor.

But the length of $S$ is greater than $4A\omega$ since the distance
along the curve through one ripple is greater than $4A$.  Thus by
choosing $A$ sufficiently small and then choosing $\omega$
sufficiently large the lemma is proved (relative to the boundary) in
this special case.

For future reference, we shall denote the 1--dimensional rippling
diffeomorphism, just constructed, by $r$.  Or to be really precise, we
use $r$ to denote this rippling diffeomorphism without phasing out
near $\{0,1\}$.

For the general case, we use induction on a handle decomposition of
$M$.  Here is a slightly inaccurate sketch of the procedure.  At each
step of the induction we move one handle keeping fixed the attaching
tube.  The handle gets replaced by a neighbourhood of its core (and
there is a compensatory enlargement of handles attached onto it)
otherwise the decomposition remains fixed throughout the induction.
For each handle in turn we think of the core as a cube, choose one
direction in the cube and one direction perpendicular to the core.  Then,
using the model in $\re^2$ just given, we ripple the chosen direction
in the core, crossing with the identity on other coordinates and
phasing out near the boundary of the core.  Then we choose another
direction in the cube and repeat this move.  This has the effect of
creating perpendicular, and smaller, ripples on top of the ripples
just made (figure 2).  We repeat this for each direction in the cube
and the end result is that all distances (in the core) are scaled.  We
then scale distances in the handle near the core by expanding a very
small neighbourhood of the core onto a small one.  Finally we redefine
the handle to lie inside the very small neighbourhood (changing the
handle decomposition of $M$) and proceed to the next handle.

\figure
\epsfxsize 3.5 truein
\epsfbox{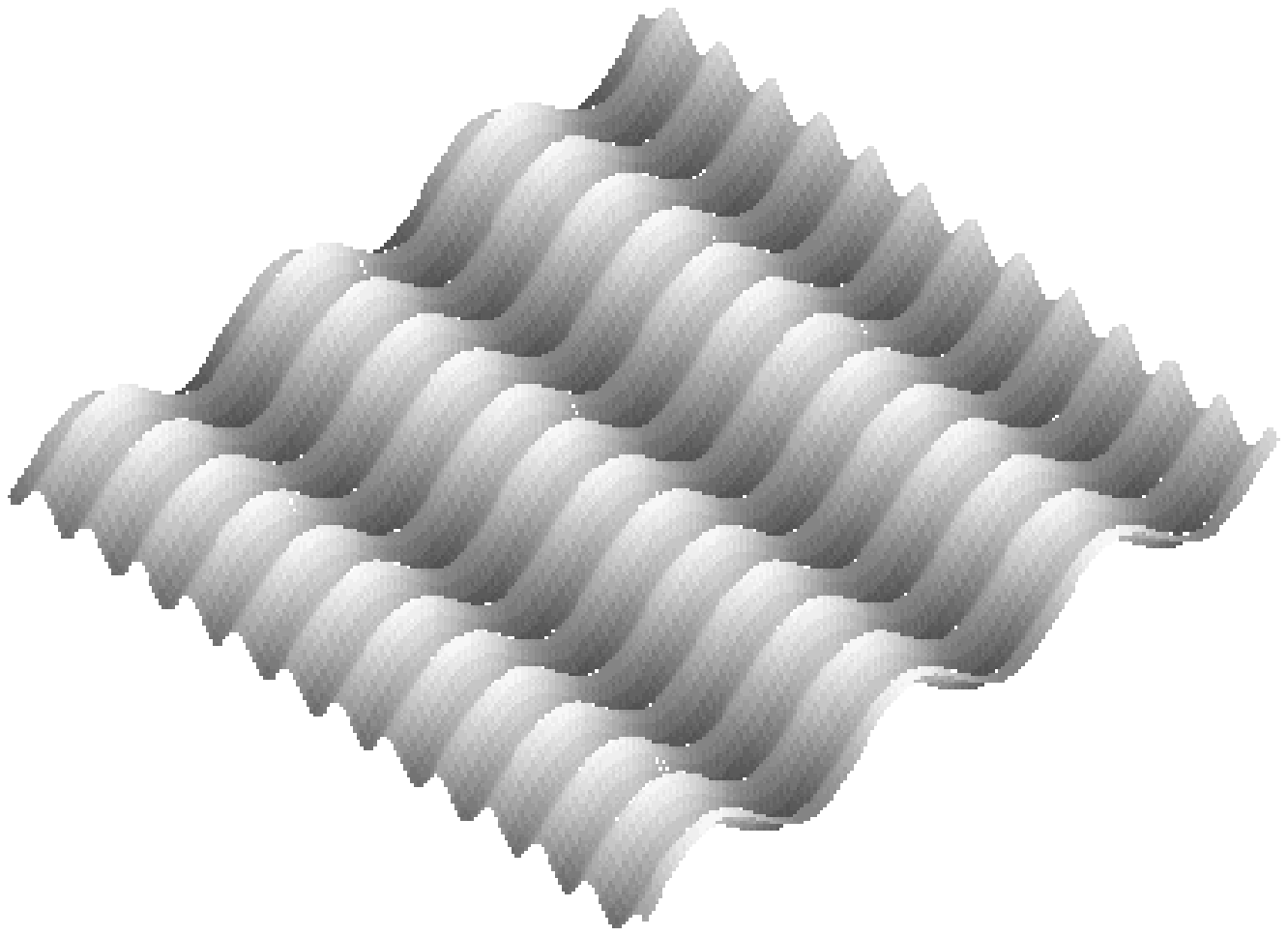}
\caption{The effect of two successive ripples in the middle of
the core}
\endfigure

\sh{The model}

For the details, let $I^j$ be the standard $j$--cube in $\re^j$.  We
shall construct a model ripple $r_j$ of $I^j$ in $I^j\times\re$ by
using the standard 1--dimensional ripple $r$, defined earlier, $j$
times.  Let $T$ be a given constant and let $I^j_0\subset {\rm
int}I^j$ be a given concentric copy of $I^j$.  Let $\phi$ be a bump
function which is 1 on $I^j_0$ and 0 outside a compact subset of ${\rm
int}I^j$.

For each $t=1,2,\ldots,j$ we perform the following inductive move.
Starting with $t=1$ we ripple $I^j$ by using $r$ on the first
coordinate and the identity on the remaining $j-1$ and phase out by
using $\phi$.  We choose the parameters for $r$ to scale distances by
the the given constant $T$.  This defines a subset $I^j_1$ of
$I^j\times\re$ and a diffeomorphism $r_1\co I^j\to I^j_1$ which
stretches distances in $I^j_0$ in the direction of the first
coordinate by the factor $T$.  Now suppose inductively that
$r_{t-1}\co I^j\to I^j_{t-1}$ has been defined.  Choose a small
orthogonal normal bundle on $I^j_{t-1}$ in $I^j\times\re$ and use
this and $r_{t-1}$ to identify a neighbourhood $U$ of $I^j_{t-1}$ with
$I^j\times V$ where $V$ is an open interval in $\re$ containing $0$.
Using this identification define $s_t\co I^j_{t-1}\to U$ by using
exactly the same construction as used for $r_1$, but replacing the
first coordinate by the $t^{\rm th}$ and choosing the height of the
sine function sufficently small that the image lies in $U$, and
adjusting the frequency so that the scale factor is again $T$.  Define
$I^j_t=s_t(I^j_{t-1})$ and $r_t=s_t\circ r_{t-1}$.

Now throughout this inductive process, the coordinate system for
$r_t(I^j)$, coming from the standard coordinate system for $I^j$,
remains perpendicular on $I^j_0$.  It follows that $r_j$ scales all
distances in $I^j_0$ by factor approximately $T$.  The reason why the
factor is not exactly $T$ is because, after the first ripple, $I^j_1$
is not flat.  Hence there is a shrinkage effect because moving out in
the $\re$--direction can move points closer together.  By choosing $U$
sufficiently small, this effect can be made as small as we please.
Hence all distances in $I^j_0$ are stretched by a factor, which may
vary with direction, but which is uniformly as close as we like to
$T$.

Another effect occurs outside $I^j_0$.  Here the non-constant scaling
in say the first coordinate (due to the phasing to zero using $\phi$)
causes the second coordinate system to become non-perpendicular and
hence the second scaling may in fact shrink some distances.  But it
can be checked that this shrinking is by a factor $\ge{\rm
sin}(\alpha)$ where $\alpha$ is the minimum angle between the images
of two lines at right angles.  But by choosing the size of all the
ripples to be small compared to the distance between $\d I^j_0$ and
$\d I^j$ (in other words the distance over which $\phi$ varies) the
distortion in the coordinate system can be made as small as we please
and hence ${\rm sin}(\alpha)$ chosen as near as we like to 1.  The way
to think of this is that we are using small ripples whose height
varies over a much larger scale.  Thus by choosing the parameters
carefully we can assume that, in the model ripple $r_j$, all distances
in $I^j_0$ are scaled up by factor as close to $T$ and as nearly
independent of direction as we please.  Outside $I^j_0$ distances are
scaled up by a factor which varies from point to point, but which, at
a given point is again as nearly independent of direction as we
please.

This completes the construction of the model.

Now to prove the lemma choose a finite (or locally finite) handle
decomposition of $M$ and suppose that $M_1=M_0\cup h_j$ (where $h_j$
is a handle of index $j$).  Suppose inductively that for given
$\delta>0$ we have constructed an isotopy finishing with a
diffeomorphism $f_0$ which moves points at most $\delta$ such that
property $(*)$ holds near $M_0$ with $f_0$ in place of $f$. 

Now choose a diffeomorphism $g$ of ${\rm int}I^j$ with the core of
$h_j$ minus attaching tube and let $I^j_0$ be chosen so $g({\rm
int}I^j-I^j_0)$ is contained in the neighbourhood of $M_0$ where $(*)$
already holds.  Choose a small orthogonal line bundle $\zeta$ on $h_j$
minus attaching tube in $Q$ and use this to extend $g$ to an embedding
of ${\rm int}I^j\times V$ where $V$ is an open interval in $\re$
containing $0$.  Define $f_1$ to be $g\circ r_j\circ g\inv$, where
$r_j$ is the model ripple contructed above.  Then, by choosing $T$
sufficiently large (noting that the directional derivatives of $g$ are
bounded on a compact subset of ${\rm int}I^j$), $f_1$ satisfies
property $(*)$ for points in a neighbourhood of $M_0$ and in the core.
There is again a shrinkage effect due to the fact that $h_j$ is not
flat (hence moving out along $\zeta$ can move points closer together).
By choosing $\zeta$ sufficiently small, this effect can again be made
as small as we please.  To stretch distances perpendicular to the core
in the handle we choose a small neighbourhood $V$ of the core
containing a much smaller neighbourhood $V'$. Then the isotopy which
expands $V'$ onto $V$ stretches distances perpendicular to the core.
Finally we change the handle decomposition of $M$ by redefining $h_j$
to lie inside $V'$ (this could be done by a diffeomorphism of $M$ and
hence defines a new decomposition).  We now have property $(*)$ in a
neighbourhood of $M_1$.  Note that the finite number of small moves
used on $h_j$ can be assumed to move points at most any given
$\delta'$.  Thus by choosing successive moves bounded by a sequence
which sums to less than $\ep$, lemma is proved by induction.
\endprf

\rk{Addenda}The proof of the ripple lemma leads at once to various
extensions: 

(1)\qua There is a codimension 0 version (ie $q=n$) as follows.  Let
$X$ be a spine of $M$.  Choose a handle decomposition with cores lying
in $X$ and no $n$--handles.  Apply the proof to this decomposition.
We obtain a neighbourhood $N$ of $X$ in $M$ and a small isotopy of $N$
in $Q$ such that $(*)$ holds in $N$.

(2)\qua In the non-compact case we can assume that $R$ and $\ep$ are
arbitrary positive functions.  (This follows at once from the local
nature of the proof.)

(3)\qua The proof gives both relative and parametrised versions.  In
the relative version we can assume that the isotopy is fixed on some
closed subset and obtain $(*)$ outside a given neighbourhood of $C$.
(This follows at once from the method of proof.) In the parametrised
version we are given a family of embeddings parametrised by a manifold
$K$ and obtain a $K$--parameter family of small isotopies such that
$(*)$ holds for each finishing embedding, where both $R$ and $\ep$ are
functions of $K$.  To see this, we use the same scheme of proof but at
the crucial stage we choose a diffeomorphism of a neighbourhood of
(the core of $h_j)\times K$ in $Q\times K$ with a neighbourhood of
$I^j\times K$ in $\re^q\times K$.  We then use the same model rippling
move over each point of $K$ varying the controlling parameters
appropriately.  The rest of the proof goes through as before.

(4)\qua Since the process is local, there is an immersed version of
the lemma in which $M$ is immersed in $Q$ and a regular homotopy is
obtained.  Further all these extensions can be combined in obvious
ways.

\rk{Proof of the flattening lemma}
We can now deduce the flattening lemma.  What we do is ripple the
horizontal coordinate using the immersed version of the ripple lemma.
This has the effect that horizontal distances are all scaled up.  We
leave the vertical coordinate fixed.  The embedding is now as flat as
we please.  The addenda to the flattening lemma follow from the
addenda to the ripple lemma given above. \qed

In order to prove the Normal Deformation Theorem (in the next section)
we shall need a bundle version of the lemma:

\proclaim{Bundle version of the flattening lemma}
Suppose that $M^n\subset W^{t+q}$ where $q-n\ge 1$ and $TW$ contains a
subbundle $\xi^t$ (thought of as vertical) such that $\xi|M$ is normal
to $M$.  Let $H$ (the horizontal subset) be the subset of $G_n(W)$ of
$n$--planes orthogonal to $\xi$.  Then given a neighbourhood $U$ of
$H$ in $G_n(W)$ there is a $C^0$--small isotopy of $M$ in $Q$ carrying
$M$ to a position where $TM\subset U$.

\prf Use a patch by patch argument.  Approximate locally as a product
$\re^t\times$patch and use the (relative) $\re^t$ version.\endprf

\rk{Remarks} (1)\qua There are obvious extensions to the bundle version 
corresponding to the extensions to the $\re^t$ version.  For the
proofs we use a similar patch by patch argument together with the
appropriate extension of the $\re^t$ version.

(2)\qua Rippling has been used on occasions by several previous
authors, in particular Kuiper used it (together with some
sophisticated estimates) to prove his version of the Nash isometric
embedding theorem [\Kuiper].

\section{Normal deformations and Gromov's theorem}

Gromov's theorem asserts (roughly) that if $M\subset W$ and the
tangential information is deformed within a neighbourhood of $M$ (ie,
$TM$ is deformed as a subbundle of $TW$) then the deformation can be
followed within a neighbourhood by an isotopy of $M$ in $W$.  We shall
prove the theorem in the following equivalent normal version (ie, $M$
follows a deformation of a bundle normal to $M$ in $W$).  The normal
version follows quickly from a repeated application of the flattening
lemma.  We shall give a precise statement (and deduction) of Gromov's
theorem after this version.

\proclaim{Normal Deformation Theorem}Suppose that $M^n\subset W^w$ and
that $\xi^t$ is a subbundle of $TW$ defined in a neighbourhood $U$ of
$M$ such that $\xi|M$ is normal to $M$ in $W$ and that $t+n<w$.
Suppose given $\ep >0$ and a homotopy of $\xi$ through subbundles of
$TW$ defined on $U$ finishing with the subbundle $\xi'$.  Then there
is an isotopy of $M$ in $W$ which moves points at most $\ep$ moving
$M$ to $M'$ so that $\xi'|M'$ is normal to $M'$ in $W$.\endproc

\prf  Suppose first that $M$ and $U$ are compact.  It follows that the
total angle that the homotopy of $\xi$ moves planes is bounded and we
can choose $r$ and a sequence of homotopies $\xi=\xi_0\simeq
\xi_1\simeq \ldots\simeq \xi_r=\xi'$ so that for each $s=1,\ldots,r-1$ the
planes of $\xi_{s-1}$ make an angle less than $\pi/4$ with those of
$\xi_s$.  We apply (the bundle version of) the flattening lemma $r+1$
times each time moving points of $M$ at most $\ep\over r+1$.  Start by
flattening $M$ to be almost perpendicular to $\xi=\xi_0$.  Then
$\xi_1$ is now normal to $M$ and we flatten $M$ to be almost
perpendicular to $\xi_1$. $\xi_2$ is now normal to $M$ and we flatten
$M$ to be almost perpendicular to $\xi_2$ etc.  After $r+1$ such moves
$M$ has been moved to $M'$ say which is almost perpendicular to
$\xi_r=\xi'$ and in particular is normal to it.

For the non-compact case, argue by induction over a countable union
of compact pieces covering $M$.\endprf

\rk{Addenda}

The proof moves $M$ to be almost perpendicular to $\xi'$ (not just
normal).  Further it can readily be modified to construct an isotopy
which ``follows'' the given homotopy of $\xi$, in other words if
$\xi_i$ is the position of $\xi$ at time $i$ in the homotopy and $M_i$
the position of $M$ at time $i$ in the isotopy, then $\xi_i|M_i$ is
normal to $M_i$ for each $i$.  To see this reparametrise the time for
the isotopy so that the flattening of $M$ with respect to $\xi_i$
takes place near time $i\over r$.  This produces a rather jerky
isotopy following the given bundle homotopy but by breaking the
homotopy into very small steps, the isotopy becomes uniform and
furthermore (apart from the initial move to become almost
perpendicular to $\xi$) $M_i$ is almost perpendicular to $\xi_i$
throughout.

There are obvious relative and parametrised versions similar to those
for the flattening lemma (and proved using those versions), and
furthermore we can assume that $\ep$ is given by an arbitrary positive
function in non-compact cases.  Finally there is a codimension 0
($t+n=w$) version which it is worth spelling out in detail, since this
is the version that implies Gromov's theorem:

\sl Suppose in the Normal Deformation Theorem that $t+n=w$ and that $M$ 
is open or has boundary and that $X$ is a spine of $M$.  Then there is
an isotopy of $M$ of the form: shrink into a neighbourhood of
$X$ followed by a small isotopy in $W$, carrying $M$ to be almost
perpendicular to $\xi'$.\rm

\proclaim{Gromov's Theorem}{\rm[\Grom; 2.4.5 $\rm(C')$, page 194]}\qua
Suppose that $M^n\subset W^w$ and suppose that $M$ is either open or
has boundary and that we are given a deformation of $TM$ over the
inclusion of $M$ to a subbundle $\eta$ of $TW$.  Then there is an
isotopy of $M$ carrying $TM$ into a given neighbourhood of $\eta$ in
the Grassmannian $G_n(W)$.\endproc

\prf Let $\xi$ be the orthogonal complement of $TM$ in $TW$.  Then
the deformation of $TM$ gives a deformation of $\xi$ to $\xi'$ say.
Pulling the bundles back over a neighbourhood of $M$ gives the
hypotheses of the Normal Deformation Theorem (codim 0 case above).
The conclusion is the required isotopy.\endprf

\rk{Final remarks}It is easy to reverse the last argument and deduce
the Normal Deformation Theorem from Gromov's Theorem.  

Another proof of the Normal Deformation Theorem is given in [\CompI] as
an extension of the arguments used to prove the Compression Theorem.
The proofs in [\CompI] are quite different in character from those
presented here.  We think of the bundles locally as defined by
independent vector fields and define flows by extending these vector
fields in an explicit fashion.  The resulting embeddings are far more
precisely defined: indeed instead of a multiplicity of ripples, there
is in the simplest case a single twist created around a certain
submanifold which we call the {\sl downset}.  For details here see the
pictures in section 3 and the arguments in section 4 of [\CompI].

We are very grateful to both Yasha Eliashberg and the referee for
comments which have greatly improved the clarity of the proof of the
ripple lemma.

\references

\bye